\documentclass[11pt]{amsart}
\usepackage{amssymb,amsmath,latexsym,amscd,amsfonts}
\usepackage{graphics}
\usepackage{epsfig}
\usepackage{setspace}

\newtheorem{thm}{Theorem}
\newtheorem{lem}[thm]{Lemma}

\newtheorem{dfn}[thm]{Definition}

\newtheorem{example}{Example}[section]

\newtheorem{conjecture}{Conjecture}

\def\hpic #1 #2 {\mbox{$\begin{array}[c]{l} \epsfig{file=#1,height=#2} \end{array}$}}
\def\vpic #1 #2 {\mbox{$\begin{array}[c]{l} \epsfig{file=#1,width=#2} \end{array}$}}
\def\ignore #1 {}
\def\cross{\times}
\def\fn{f^{[n]}}

\def\Z{\mbox{{$\mathbb{Z}$}}}

\def\E{\mbox{$\exists$}}
\def\A{\mbox{$\forall$}}

\def\B{\mbox{${\mathcal B}$}}

\def\ss{\subset}
\def\pf{\mbox{Proof: }}
\def\de{\delta}
\def\th{\theta}
\def\a{\alpha}

\def\qed{\mbox{$\Box$}\\}

\newcommand{\etilde}{\widetilde{e}}
\newcommand{\ftilde}{\widetilde{f}}
\newcommand{\mutilde}{\widetilde{\mu}}
\renewcommand{\E}{{\mathbb E}}
\newcommand{\R}{{\mathbb R}}
\renewcommand{\d}{\delta}
\newcommand{\g}{\gamma}

\newcommand{\x}{{\mathbf x}}
\newcommand{\y}{{\mathbf y}}

\renewcommand{\A}{{\mathcal A}}
\newcommand{\e}{\epsilon}
\def\Tfn{T_{\mu,n}}
\def\Sfn{S_{\mu,n}}

\def \aaron #1 {\marginpar{#1 -AA}}

\begin{document}


\title{Optimal Estimators for Threshold-Based Quality Measures\footnote{This article was published in \emph{Journal of Probability and Statistics,} vol.~2010, Article ID 752750 (2010).}}
\author{Aaron Abrams, Sandy Ganzell, Henry Landau, Zeph Landau, James Pommersheim, 
and Eric Zaslow}
\begin{abstract}
We consider a problem in parametric estimation:  given $n$ samples from an
unknown distribution, we want to estimate which distribution, from a given 
one-parameter family, produced the data.  Following Schulman and Vazirani
\cite{schulman}, we evaluate an estimator in terms
of the chance of being within a specified tolerance of the correct answer, in
the worst case.  We provide optimal estimators for several families of distributions 
on $\R$.  We prove that for distributions on a compact space, there is always an
optimal estimator that is translation-invariant, and we conjecture that this conclusion 
also holds for any distribution on $\R$.  By contrast, we give an example showing it 
does not hold for a certain distribution on an infinite tree.
\end{abstract}

\maketitle

\pagestyle{myheadings}
\markboth{Abrams, Ganzell, Landau, Landau, Pommersheim, and Zaslow}
{Optimal Estimators for Threshold-Based Quality Measures}

\section{Introduction}

Estimating probability distribution functions is a central problem in statistics.
Specifically, beginning with an unknown probability distribution on an underlying
space $X$, one wants to be able to do two things:  
first, given some empirical data sampled from the unknown 
probability distribution, estimate which one of a presumed set of 
possible distributions produced the data; and second, obtain 
bounds on how good this estimate is.  
For example, the \emph{maximum likelihood estimator} selects the distribution
that maximizes the probability (among those under consideration) of producing 
the observed data.  Depending
on what properties of the estimator one is trying to evaluate, this may or may 
not be optimal.  An extensive literature, dating back to the early 20th century,
addresses problems of this sort; see for example 
\cite{bondar,eaton,fisher,kiefer,pitman,wijsman}.

%

In this paper we consider one such problem.  We presume samples are coming 
from an unknown ``translate" of a fixed known distribution.  The challenge is to 
guess the translation parameter.  More precisely, we are given a distribution 
$\mu$ on a space $X$, along with an action of a  group $G$ on $X$, which defines 
a set of translated distributions $\mu_{\th}$ as follows:
\begin{equation}\label{Def:Shift}
\mu_{\th}(A)=\mu(\{x:\th x\in A\})
\end{equation}
for $A\subset X$.
%
Thus in this context an \emph{estimator} is a (measurable) function $e:X^n \to G$; the input
$\x=(x_1,\ldots,x_n)$ is the list of samples, and the output $e(\x)$ 
is the estimate of $\th$, the translation parameter.  For the majority of the paper we will
study the case of $G=\R$ acting by translations (changes in location) on $X=\R$, and the group action will be written additively, as seen beginning in Section \ref{Sec:SBound}.

%

We are interested in finding good estimators; thus we need a way 
of measuring an estimator's quality.  A common way to do this 
is to measure the \emph{mean squared error}, in which case an
optimal estimator minimizes this error.  Various 
results are known in this case; for instance the maximum likelihood
estimator (which agrees with the \emph{sample mean estimator}),
$$e(\x)=\left( \frac 1 n \sum x_i \right) - \E(\mu),$$
minimizes the mean squared error if $\mu$ is a Gaussian distribution on $\R$.

In this paper we investigate a different and natural measure of quality whereby
we consider an estimator to succeed or fail according to whether or not its
estimate is within a certain threshold $\d>0$ of the correct answer.  We then define the 
quality of the estimator to be the chance of success in the worst case.
This notion was introduced in \cite{schulman} to analyze certain approximation algorithms 
in computer science.
Precisely, the \emph{$\d$-quality} of $e$ is defined as
\begin{eqnarray}\label{Def:Quality}
Q_{\d}(e)&=& \inf\limits_{\th} Q_{\d}^{\th} (e)\\
&=&
\inf\limits_{\th} \Pr \left\{ d(e(\x),\th) < \d : x_i\mbox{ are chosen from } \mu_{\th} \right\}  
\nonumber\\
&=&
\inf\limits_{\th} \mu_{\th}^n(\{\x : d(e(\x),\th) < \d \}), \nonumber
\end{eqnarray}
where $d$ is a metric on $X$ and $\mu_{\th}^n$ is the product measure 
$\mu_{\th}\cross\cdots\cross\mu_{\th}$ on $X^n$.\footnote{In
the case of perverse measures, $\mu,$ we must consider the probability
as the sup of the intersection of the set $\{d(e(\x),\th)<\delta\}$ with all measurable
sets.  We will ignore this caveat throughout. Indeed, we primarily focus on
absolutely continuous measures (as \cite{farrell} and \cite{kamberova} have done, 
for example) and purely atomic measures.}  Note that, depending on context, it
is sometimes advantageous to define the quality using a closed interval rather
than an open one; for example in the discrete case we could then interpret 
$Q^{\th}_0(e)$ as the probability that $e$ is exactly equal to $\th$.
We write $Q(e)$ when the
value of $\delta$ is unambiguous.  For fixed $\d$, an ($n$-sample) estimator 
$e$ is \emph{optimal} if $Q_\d(e)\geq Q_\d(e')$ 
for all ($n$-sample) estimators $e'$.  Many authors use the term \emph{minimax}
to describe optimal estimators.  Note that much of the literature on this
subject uses the notion of loss functions and the associated risk $R=1-Q$;
our point of view is equivalent but more convenient for our purposes.

Motivated initially by analyzing an approximate
algorithm for determining the average matching size in a graph, Schulman and Vazirani 
\cite{schulman} introduce the stronger notion of a \emph{majorizing estimator},
which is optimal (by the above definition) simultaneously for all $\d>0$.  This was 
previously studied by Pitman \cite{pitman}, who considered several different optimality criteria 
and, for each one, constructed optimal ``shift-invariant'' estimators (defined below).
Schulman and Vazirani focus on the Gaussian distribution and prove that the mean
estimator is the unique majorizing estimator in this case. 



In the first part of this paper we investigate the optimal estimators for several different 
classes of distributions on $\R$.  We conjecture that there is always
an optimal estimator $e$ that is \emph{shift-invariant}, i.e.~$e$ satisfies
$$e(x_1+c,\ldots,x_n+c)=e(x_1,\ldots,x_n)+c$$
for all $c,x_i\in\R$.  These estimators are typically easier to analyze than
general estimators, because the quality is the same everywhere, 
i.e.~$Q(e)=Q^\th(e)$ for every $\th$.
Conditions under which invariant minimax estimators can be obtained have been 
studied, for example, in \cite{berger}, \cite{lehman} and \cite{robert}.  Indeed, some 
of our existence results follow from the quite general Hunt-Stein theorem 
\cite[Theorem 9.5.5]{robert}, but we give constructions that are very natural and explicit.
We obtain general bounds on the quality of shift-invariant estimators
(Section \ref{Sec:SBound}) and general estimators (Section \ref{Sec:TBound}), 
and then we apply these bounds to several families of distributions 
(Section \ref{Sec:Examples}).  In each case, we are 
able to construct an optimal estimator that is shift-invariant.  These 
examples include the Gaussian and exponential distributions, among others.

These results motivate our study of shift-invariant estimators on other spaces;
these are estimators that are equivariant with respect to the induced
diagonal action of $G$ on either the left or the right on $X^n$.  That is,
a \emph{left-invariant} estimator satisfies
\begin{equation}\label{Def:ShiftInv}
e(g\x)=ge(\x)
\end{equation}
where 
$$g(x_1,\ldots,x_n)=(gx_1,\ldots,gx_n).$$
Right-invariance is defined similarly.

In Section \ref{Sec:Compact} we show that on a compact space $X$, if $e$ 
is an estimator for $\mu$, then there is always a
shift-invariant estimator with quality at least as high as that of $e$.
The idea is to construct a shift-invariant estimator $s$ as an average
of the translates of $e$; this is essentially a simple proof of a special case
of the Hunt-Stein theorem.  As there is no invariant probability measure 
on $\R$, the proof does not extend to the real case.

Finally, in the last section, we give an example due to L.~Schulman
which shows that (on non-compact spaces) there may be no shift-invariant
estimator that is optimal.  It continues to be an interesting problem to 
determine conditions under which one can guarantee the existence of 
a shift-invariant estimator that is optimal.

The authors thank V.~Vazirani and L.~Schulman for suggesting the problem 
that motivated this paper along with subsequent helpful discussions. We are 
grateful to the reviewers for many helpful suggestions, especially for correcting 
our proof of Lemma \ref{Lem:Decompose}.  As always, 
we thank Julie Landau, without whose support and down-the-line backhand this 
work would not have been possible.

\section{The real case:  shift-invariant estimators}
\label{Sec:SBound}

Let $G=X=\R$, and consider the action of $G$ on $X$ by translations.
Because much of this paper is concerned with this context, we spell out
once more the parameters of the problem.  We assume $\d>0$ is fixed 
throughout.  We are given a probability
distribution $\mu$ on $\R$, and we are to guess which distribution $\mu_\th$
produced a given collection $\x=(x_1,\ldots,x_n)$ of data, where 
$\mu_\th(A)=\mu(\{x : x+\th\in A\})$.
An estimator is a function $e:\R^n\to\R$, and we want to maximize its quality,
which is given by
\begin{eqnarray*}
Q(e)=\inf\limits_{\th} Q^{\th} (e)&=&
\inf\limits_{\th} \Pr \left\{ |e(\x)-\th| < \d : x_i
\mbox{ are chosen from } \mu_{\th} \right\} \\
&=&\inf\limits_{\th} \mu_{\th}^n(\{\x : |e(\x)-\th| < \d \}).
\end{eqnarray*}

First some notation.  We will write the group action additively and
likewise the induced diagonal action of $G$ on $\R^n$;
in other words if $\x=(x_1,\ldots,x_n) \in \R ^n$ and $a \in \R$, then 
$\x+a$ denotes the point $(x_1+a,\ldots,x_n+a)\in\R^n$.  Similarly
if $Y\subset \R^n$ and $A\subset \R$ then $Y+A=\{\y+a : \y\in Y, a\in A\}$.
We also use the ``interval notation'' $(\x+a,\x+b)$ for the set $\{\x+t : a<t<b\}$;
this is a segment of length $(b-a)\sqrt{n}$ in $\R^n$ if $a$ and $b$ are finite.
If $f: \R^n \rightarrow \R$ is any function, and $\theta \in \R$, define
$f_{\th}(\x)=f(\x-\th)$.  
If $f:\R\to\R$, then define $f^{[n]}:\R^n\to\R$ by 
$f^{[n]}(\x)=f(x_1)f(x_2)\cdots f(x_n)$.

We now establish our upper bound on the quality of shift-invariant estimators.
Note that a shift-invariant estimator has the property that $e(\x -e(\x))=0$.
Also note that a shift-invariant estimator is determined uniquely by its values
on the coordinate hyperplane 
$$X_0=\{\x\in \R^n : x_1=0 \},$$ 
and that a shift-invariant estimator exists for any choice of such values on $X_0$.
In addition, for $e$ shift-invariant,
$$\mu_{\th}^n \left( \{\x: |e(\x)-\theta | < \delta \}\right) =
\mu_{\th}^n \left( \{\x: |e(\x-\theta)| < \delta \}\right) =
\mu^n \left( \{\x: |e(\x)| < \delta \}\right),$$
so the quality can be ascertained by setting $\theta = 0.$

\begin{dfn} \label{A}
For fixed $n$, let ${\mathcal A}$ denote the collection of all 
Borel subsets $A$ of the form
\[
A=\bigcup_{\x\in X_0}\x+\left(f(\x)-\d,f(\x)+\d\right),
\]
where $f:X_0\rightarrow\R$ is a Borel function.
For fixed $\mu$ and $n$, define
$$\Sfn = \Sfn(\d)=\sup_{A \in {\mathcal A}} \{\mu^n(A)\}.$$
\end{dfn}                                    

\begin{thm}
\label{sbound}
Let $\mu$ and $n$ be given.  Then any shift-invariant $n$-sample estimator $e$ 
satisfies $Q(e)\leq \Sfn$.
\end{thm}
 
\pf  Due to the observation above, it suffices to bound the quality of 
$e$ at $\th=0$.  But this quality is just $\mu^n(A)$ where 
$A=e^{-1}((-\d,\d))$.  Note that 
$$A = \bigcup\limits_{\x \in X_0} (\x-e(\x) - \de, \x-e(\x) +\de), $$
and in particular $A\in\A$.
Thus the quality of $e$ is at most $\Sfn$.
\qed

\begin{thm}
\label{shifty}
Let $\mu$ and $n$ be given.  If the $\sup$ in Definition~\ref{A} is 
achieved, then there is a shift-invariant $n$-sample estimator with
quality $\Sfn$.  For any $\e>0$, there is a shift-invariant $n$-sample 
estimator $e$ with quality greater than $\Sfn-\e$.
\end{thm}

\pf For a given $A \in {\mathcal A}$, let $f$ be the corresponding Borel
function (see Definition \ref{A}).  Define the estimator $e$ to be $-f(\x)$
on $X_0$, and then extend to all of $\R^n$ to make it shift-invariant.
Note that $A \ss e^{-1}((-\d,\d))$, so $Q(e)\geq \mu^n(A)$.  The 
theorem now follows from the definition of $\Sfn$.
\qed

\section{The real case:  general estimators}
\label{Sec:TBound}

In this section we obtain a general upper bound on the quality of
randomized estimators, still in the case $G=X=\R$.  The arguments are 
similar to those of the previous section.

Again $\d$ is fixed throughout.  A \emph{randomized estimator}
is a function $X^n\times\Omega\to\R$ where $\Omega$ is a probability
space of estimators; thus for fixed $\omega\in\Omega$, $e=\etilde(\cdot,\omega)$
is an estimator.
The $\d$-quality of a 
randomized estimator $\etilde$ is
$$Q(\etilde)=\inf\limits_\th Q^\th (\etilde)$$
where
$$Q^\th (\etilde) = \E \mu_\th^n \left\{\x \mid \left|\etilde(\x)-\th\right|<\d \right\}.$$

\begin{dfn} \label{B}
For fixed $n$, let 
$$\B = \{B \ss \R^n : B \cap (B+ 2k\d) = \emptyset \ 
\text{for all nonzero integers } k \}.$$
For fixed $\mu$ and $n$, define
$$\Tfn =\Tfn(\d) = \sup_{B \in \B} \{ \mu^n(B) \}.$$
\end{dfn}

Comparing with Definition \ref{A}, we observe that $\A\subset \B$ and 
hence $\Sfn \leq \Tfn$.

\begin{thm}\label{upbd}
Let $\mu$ and $n$ be given.  Then any $n$-sample
randomized estimator $\etilde$ satisfies
$Q(\etilde) \leq \Tfn$.
\end{thm}

\pf
We will give a complete proof in the case that $\mu$ is defined by a density
function $f$, and then indicate the modifications required for the general case.
The difference is purely technical; the ideas are the same.

Consider first a non-randomized estimator $e$.  
The performance of $e$ at 
$\th$ is $\mu_{\th}^n \left( e^{-1}((\th - \de, \th + \de)) \right)$.  
To simplify notation we will let $I_{e,\th}$ denote the set 
$e^{-1}((\th - \de, \th + \de))$ and we will suppress the subscript 
$e$ when no ambiguity exists.  Since $Q(e)$ is an infimum,
the average performance of $e$ at the $k$ points $\th_i= 2 \de i$, 
($i=1,2, \dots, k$) is at least $Q(e)$:
\begin{equation} \label{e1}
Q(e) \leq \frac{1}{k} \sum _{i=1}^k \mu_{\th_i}^n( I_{\th_i})
\end{equation}

Now we use the density function $f$.  Recall that $f_{\th}(x)=f(x-\th)$.
Define $\ftilde$ on $\R^n$ by
$$\ftilde(\x)= \max _i \{\fn_{\th_i}(\x)\}= \max_i
\{f_{\th_i}(x_1) f_{\th_i}(x_2) \cdots
f_{\th_i}(x_n) \}.$$  
%
%
%
%

Since the $I_{\th_i}$ are disjoint, we now have
\begin{eqnarray}
\label{e2}
Q(e) &\leq& 
\frac{1}{k} \sum_{i=1}^k \int_{I_{\th_i}} \fn_{\th_i}(\x) d\x \nonumber\\
 &\leq&\frac{1}{k}\sum_{i=1}^k \int_{I_{\th_i}}\ftilde(\x)d\x \nonumber\\
 &\leq& \frac 1 k \int_{\bigcup I_{\th_i}} \ftilde(\x) d\x \nonumber\\
 &\leq& \frac 1 k \int_{\R^n} \ftilde(\x) d\x \nonumber\\
&=& \frac 1 k \int\limits_{\{x_1\leq0\}} 
\ftilde(\x) d\x 
+ \frac 1 k \int\limits_{\{0<x_1<2\d k\}} {\ftilde(\x) d\x}
+ \frac 1 k \int\limits_{\{x_1\geq2\d k\}}
\ftilde(\x)d\x.
\end{eqnarray}

We will bound the middle term by $\Tfn$ and show that the first 
and last terms go to zero (independently of $e$) as $k$ gets large.
The bound on the middle term is a consequence of the following claim.\\

{\bf Claim.}  For any $a\in\R$, $$\int_{\{a \leq x_1 \leq a+2\de\}}
\ftilde(\x)d\x \leq \Tfn.$$

To prove the claim, set $Z=\{ \x\in\R^n : a\leq x_1 < a+2\d \}$, and set 
$Z_i=\{ \x\in Z : i \text{ is the smallest index such that } \ftilde(\x)= \fn_{\th _i}(\x) \}$.  
Thus the $Z_i$ are disjoint 
and cover 
$Z$.  Now
\begin{eqnarray*}
\int_{\{a\leq x_1 \leq a+2\de\}} \ftilde (\x)d\x &=& \sum_i \int_{Z_i} \ftilde(\x) d\x 
=\sum_i \int_{Z_i} \fn_{\th_i}(\x) d\x \\
&=& \sum_i \int_{Z_i - \th_i} \fn(\x) d\x \\
&=& \int_{\bigcup (Z_i - \th_i)}\fn(\x) d\x \\
&\leq& \Tfn . 
\end{eqnarray*}
The last equality follows from the fact that the $Z_i - \th_i$ are disjoint 
(recall that $\th_i=2\d i$), and the final step follows because the set 
$B= \cup (Z_i - \th_i)$ is in $\B$.  This proves the claim. 
\qed

Next we show that $\frac{1}{k}  \int_{\{x_1\leq 0\}} \ftilde(\x)d\x$ 
approaches zero as $k$ grows.  Recall that $\th_i=2\d i$, and set $z_i= 
\int_{\{x_1\leq 0\}} \fn_{2\d i} (\x) d\x = 
\int_{\{x_1 \leq -2\d i\}} \fn(\x)d\x$.  The function $f$ is a probability density function,
so $f$ is nonnegative and has total integral $1$.
The Dominated Convergence Theorem then implies that the sequence $\{z_i\}$ is decreasing to 0.
Bounding $\ftilde (\x)$ by $\sum_{i=1}^{k} \fn_{\th_i}(\x)$ we have
\begin{equation*} \label{e3}
\frac{1}{k} \int_{\{x_1 \leq 0\}} \ftilde(\x) d\x \leq 
\frac{1}{k} \int_{\{x_1 \leq 0\}} \sum_{i=1}^{k} \fn_{\th_i}(\x)d\x =
\frac{1}{k} \sum_{i=1}^k z_i 
\to 0.
\end{equation*}

A similar argument shows that the term 
$\frac{1}{k} \int_{\{x_1 \geq 2\d k\}} \ftilde(\x)d\x$ goes to 0 as $k$ 
grows.  Since (\ref{e2}) holds for all $k$, we have $Q(e)\leq \Tfn$
for any estimator $e$.

We have shown that for any $\e>0$, we can find $k$ depending on $\e$ and 
$f$ such that the average 
performance of an arbitrary estimator $e$ on the $k$ points $\th_i=2\d i$
is bounded above by $\Tfn+2\e$.  Now, for a randomized estimator $\etilde$,
the quality is bounded above by its average performance on the same $k$ 
points, and that performance can be no better than the best estimator's 
performance.  We conclude that $Q(\etilde) \leq \Tfn + 2 \e$,
and the theorem follows.  

The proof is now complete in the case that $\mu$ has a density $f$.  In general,
the argument requires minor technical adjustments.  The first step that requires
modification is the definition of the function $\ftilde$.
Let
\[
\nu =\sum^{k}_{i=1}\mu_{\theta_i}^{n},\quad g_i=\frac{d\mu_{\theta_i}^{n}}{d\mutilde}\quad\mbox{and}\quad \ftilde=\max(g_1,\dots g_k).
\]
Then $\ftilde\cdot\nu=\mutilde$ and we work with $\mutilde$ rather than $\ftilde$.
From here one defines the $Z_i$ accordingly, and the remainder of the
argument goes through with corresponding changes.
\qed

\section{The real case:  examples}
\label{Sec:Examples}

We have obtained bounds on quality for general estimators and for 
shift-invariant ones.  In this section we give several situations where
the bounds coincide, and therefore the optimal shift-invariant 
estimators constructed in Section~\ref{Sec:SBound} are in fact optimal estimators,
as promised by the Hunt-Stein theorem.
These examples include many familiar distributions, and they provide evidence 
for the following conjecture.

\begin{conjecture}
Let $\mu$ be a distribution on $\R$.  Then there is an optimal estimator
for $\mu$ that is shift-invariant.
\end{conjecture}

\subsection{Warm-up:  unimodal densities, one sample}

Our first class of examples generalizes Gaussian distributions and many others.
The argument works only with one sample, but we will refine it in \ref{Subsec:Sufficient}.  
Note that the optimal estimator in this case is the maximum likelihood estimator.

We say that a density function is \emph{unimodal} if for all $y$, $\{x:f(x)\geq y\}$ is
convex.

\begin{example}\label{unimodal}
Let $\mu$ be defined by a unimodal density function $f$.  
Then there is a shift-invariant one-sample estimator that is optimal.
\end{example}

\pf  We first show that $T_{\mu,1}=S_{\mu,1}$.  It follows from the definition 
of ${\mathcal B}$ that any set 
$B \in {\mathcal B}$ must have Lebesgue measure less than or equal to $2\de$.  
Since $f$ is unimodal, $\int_B f(x)dx$ is maximized by concentrating $B$ around the 
peak of $f$; thus the best $B$ will be an 
interval that includes the peak of $f$.  But any interval in ${\mathcal B }$ 
is contained in ${ \mathcal A }$ and thus $T_{\mu,1}\leq S_{\mu,1}$.  
Since $S_{\mu,1} \leq T_{\mu,1}$ always, we have $T_{\mu,1}=S_{\mu,1}$. 

Now, recalling that $S_{\mu,1}$ and $T_{\mu,1}$ are defined as suprema,
we observe that the above argument shows that if one is achieved then so is
the other.  Therefore the result follows from Theorems~\ref{shifty} and \ref{upbd}.
\qed

\subsection{A sufficient condition}\label{Subsec:Sufficient}

The next class is more restrictive than the preceding, but with the
stronger hypothesis we get a result for arbitrary $n$.  Any 
Gaussian distribution continues to satisfy the hypothesis.

\begin{example}
Let $\mu$ be a distribution defined by a density function of the form $f=e^{\lambda (x)}$ 
with $\lambda ' (x)$ continuous and decreasing.  Then for any $n$, there is a shift-invariant 
$n$-sample estimator that is optimal.
\end{example}

\pf
For any fixed $\x \in X_0$, we define a function $h_\x:\R\longrightarrow\R$ by
$$
h_\x(t) = \fn(\x+t) = e^{\lambda(x_1+t)+\cdots+\lambda(x_n+t)}.
$$
Since
$$
h_\x'(t) = e^{\lambda(x_1+t)+\cdots+\lambda(x_n+t)}
(\lambda'(x_1+t)+\cdots+\lambda'(x_n+t))
$$
and $\lambda'$ is decreasing, it is clear that for each $\x$,
$h_\x'(t) =0$ for at most one value of $t$. Since
$h_\x(t) \to 0$ as $t \to \pm\infty$, it
follows that for any $\x$, $h_\x$ is a unimodal
function of $t$.
 
Now the argument is similar to Example \ref{unimodal}.
We will show that $\Tfn=\Sfn.$  Since
$\fn$ restricted to each orbit $\x+\R$ is unimodal as we have just shown,
a set $B\in {\mathcal B}$ on which the integral of $\fn$ is maximized
is obtained by choosing an interval from each orbit.  To
make this more precise, for each $\x\in X_0$, let $t_{\x}$ be the center
of the length $2\delta$ interval $I=(t_{\x}-\d,t_{\x}+\d)$ that
maximizes $\int_I h_\x\,dt$.  Then let
$$ A =\bigcup_{\x\in X_0} (\x + t_x-\delta,\x+t_x+\delta).  $$
Now $A\in{\mathcal A}$, and moreover $\mu^n(A) \geq \mu^n(B)$ for any 
$B \in \mathcal B$, because $\int_{A\cap (\x+\R)} \fn \geq
\int_{B\cap (\x+\R)} \fn$ for each $\x\in X_0$.
 
Thus $\sup_{B\in\mathcal B} \{\mu^n(B)\}$ is achieved by $B=A\in{\mathcal A}$, 
and it follows that $\Sfn=\Tfn$ and that
the best shift-invariant estimator is optimal.
\qed

\subsection{Monotonic distributions on $\R^+$}

The third class of examples generalizes the exponential distribution, defined by the 
density $f(x)=\lambda e^{-\lambda x}$ for $x \geq 0$ and $f(x)=0$ for $x<0$.
The optimal estimator in this case\footnote{Note that in a typical estimation problem 
involving a family of exponential distributions, one is trying to estimate the ``scale'' 
parameter $\lambda$ rather than the ``location'' $\th$.}
is \emph{not} the maximum likelihood estimator.

\begin{example}\label{slide}
Let $\mu$ be defined by a density function $f$ that is decreasing for $x\geq 0$
and identically zero for $x < 0$.  Then for any $n$, there is a 
shift-invariant $n$-sample estimator that is optimal.
\end{example}

\pf
We construct the estimator as follows:  for $\x \in \R^n,$ 
define $e(\x) = \min \{x_1,...,x_n\} - \delta.$
Note that this is shift-invariant; therefore $Q(e)$ can be 
computed at $\theta = 0.$
That is, it suffices to show that $Q^0(e) = \Tfn.$

Let $B = \{ \x \in \R^n : 0 \leq \min \{x_1,...,x_n\} < 2\delta\}.$  Note 
that $B = e^{-1}([-\delta,\delta)),$ and so $\mu^n(B)$ is the
quality of $e$.
Note also that $B \in {\mathcal B}$
(in fact $B\in{\mathcal A}$), so certainly $\mu^n(B)\leq \Tfn$.
We will show that any $C\in{\mathcal B}$ can be modified
to a set $C'\in{\mathcal B}$ such that $C'\subset B$ and 
$\mu^n(C)\leq \mu^n(C')$.
It then follows that $\Tfn\leq \mu^n(B)$, and this will
complete the proof.

So, let $C\in{\mathcal B}$, and define 
$C' = \{ \x \in B : \x + 2k\delta \in C \hbox{ for some } k \in \Z\}.$
Note that $k$ is determined uniquely by $\x$.
Now $C' \subset B$ is in ${\mathcal B}$, and
by our hypotheses on $f$, if $\x\in B$ then $\fn(\x)\geq \fn(\x+2k\d)$ for
every integer $k$.  Therefore $\mu^n(C') - \mu^n(C) = 
\int_{C'}[\fn(\x) - \fn(\x + 2k\delta)] d\mu^n \geq 0$.
\qed

\subsection{Discrete distributions}

Here we discuss purely atomic distributions on finite sets of points.  Because
we are only trying to guess within $\d$ of the correct value of $\th$, there are
many possible choices of estimators with the same quality.  Among the optimal
ones is the maximum likelihood estimator.

\begin{example}\label{discrete}
Let $\mu$ be a distribution on a finite set of points $Z$.  There is a shift-invariant
one-sample estimator that is optimal.  Furthermore, if all of 
the pairwise distances between points of $Z$ are distinct, then for every $n$ 
there is a shift-invariant $n$-sample estimator that is optimal.
\end{example}

\pf
We first treat the case $n=1$. Since 
$\mu$ is discrete, the supremum defining $S_{\mu,1}$ is
attained; therefore by Theorems~\ref{shifty} and \ref{upbd}
it suffices to show that every estimator has quality at most $S_{\mu,1}$.

Let $Z = \{z_1,\ldots,z_r\}$, and for any $z\in Z$, let $p_z$ denote the mass at $z$.  
For a finite set, we use $|\cdot |$ to denote the cardinality.
Suppose that $e$ is any estimator.

\begin{lem}
Let $Y$ be any finite subset of $\R$.  Then 
$$Q(e) \leq S_{\mu,1} \frac{|Y+Z|}{|Y|},$$
where $Y+Z$ denotes the set $\{y+z \mid y\in Y, z\in Z\}$.
\end{lem}

{Proof of Lemma:}  
To prove the lemma we estimate the average quality $Q^{\theta}(e)$ over $\th\in Y$. 
We have
$$\sum_{\th \in Y} Q^{\th}(e) = \sum_{\th \in Y} \mu(e^{-1}(\th-\d, \th +\d) - \th) =
\sum_{\th \in Y} \sum_{z} p_z
$$ 
with the inner part of the last sum taken over those $z\in Z$ that lie in 
$e^{-1}(\th-\de, \th +\de) -\th$.   Using $x$ to denote $\th + z$,  this condition
becomes $e(x) \in (\th-\d,\th+\d)$, and  the right hand side above may be rewritten
as 
$$\sum_{\th \in Y} \sum_{z} p_z=\sum_{\th\in Y} \sum_x p_{x-\th}
=\sum_{x \in Y+Z} \sum_{\th} p_{x-\th},$$
with the inner sum now taken over all $\th$ with
$e(x) \in (\th-\d, \th+\d)$.  
This latter condition implies that $z$ is within $\de$ of $x-e(x)$.
But by definition, $S_{\mu,1}$ is the
maximum measure of any interval of length $2\de$.  Hence, for any fixed $x\in Y+Z$, 
the inner sum is at most $S_{\mu,1}$, and the entire sum is
thus bounded above by $S_{\mu,1} \cdot |Y+Z|$.  Dividing by $|Y|$ gives a bound for
the average quality over $Y$, and since $Q(e)$ is defined as an infimum the 
lemma follows.
\qed

We now apply the lemma to complete the Example.  Let $k$ be a natural number, and let
$$Y_k = \{ h_1 z_1 + \cdots + h_r z_r : h_i\in\Z \text{ and } 0 \leq h_i < k \}.
$$
Note that $Y_k \subseteq Y_{k+1}$ and $|Y_k| \leq k^r$.  It follows that for any
$\epsilon > 0$ there exists $k$ such that $|Y_{k+1}|/|Y_k| < 1 + \epsilon$,
for otherwise $|Y_k|$ would grow at least exponentially in $k$.  Using the fact that
$Y_k+Z \subseteq Y_{k+1}$, the lemma applied to $Y_k$ implies that 
$Q(e) \leq S_{\mu,1} (1+\epsilon)$.  Therefore $Q(e)\leq S_{\mu,1}$, and
this finishes the case $n=1$.

Lastly, we consider an arbitrary $n$.  If we are given samples $x_1, ..., x_n$ and if
any $x_i \neq x_j$ for some $i$ and $j$, then by our hypothesis the shift $\th$ is uniquely
determined.  Thus we may assume that any optimal estimator picks the right
$\th$ in these cases, and the only question is what value the estimator 
returns if all the samples are identical. The above analysis of the one
sample case can be used to show that any optimal shift-invariant estimator is optimal.
\qed


\section{The compact case}
\label{Sec:Compact}


So far we have dealt only with distributions on $\R$, where the shift parameter
is a translation.  In every specific
case that we have analyzed, we have found a shift-invariant estimator among
the optimal estimators.  
%
%
%
%
%
In this section we prove that if $G=X$ is a compact group acting on itself
by (left) multiplication, then at least for measures defined by density functions,
there is always a shift-invariant estimator as good as any given estimator.
In Section~\ref{Sec:Counterexample} we show that the compactness 
hypothesis cannot be eliminated entirely; we do not know how much
it can be weakened, if at all.

We will continue to use both $G$ and $X$ as notations, in order to 
emphasize the distinction between the two roles played by this object. 
Eaton \cite{eaton} discusses estimators in a context in which the group $G$ acts 
on both the parameter space $\Theta$ and the sample space $X$.  In his work, the 
sample space $X$ is an arbitrary homogeneous space (i.e., a space with a transitive 
$G$-action).  In this generality, shift-invariant estimators may not exist, since there may
not even exist a function from $X^n$ to $\Theta$ that preserves the $G$ action. For this
reason, we choose to identify the sample space with the group $G$.

As usual $G$ acts diagonally on $X^n$; we denote
the orbit space by $Y$.  An element $\y$ of $Y$ is an equivalence
class $\y=[\x]=\{(gx_1,\ldots,gx_n):g\in G\}$, which we identify
with $G$ via $(gx_1,\ldots,gx_n) \mapsto gx_1$.
For $\x=(x_1,\ldots,x_n)\in X$ we denote by $\x_0$ the point
$x_1^{-1}\x$; thus $\x_0$ is in the orbit of $\x$ and has first
coordinate $1$.  The set $X_0=\{\x_0 : \x\in X\}\ss X^n$ is naturally
identified with $Y$.

Equip $G=X$ with a left- and right-invariant metric $d$, meaning that
$d(gx,gy)=d(x,y)=d(xg,yg)$ for all $x,y,g\in G$. 
Let $B(g)=B_\d(g)$ denote the ball of radius $\d$ around $g\in G$.
If $S$ is a subset of a measure space $(T,\alpha)$ then we denote the
measure of $S$ variously by $\alpha(S)$, $\int_S \,d\a$, or $\int_T \chi_S \,d\a$.
(The notation $\chi_S$ refers to the characteristic function of the set $S$.)

Fix $\d$ and $n$, and let $\mu$ be an arbitrary measure on $X$.
The following technical lemma says that to 
evaluate an integral over $X^n$, we can integrate over each $G$-orbit and 
then integrate the result over the orbit space.

\begin{lem}\label{Lem:Decompose}
There exist measures $\nu$ on $Y$ and
$\a_\y$ on each orbit $\y$ such that for any function $F$ on $X^n$,
$$\int_{X^n} F(\x) \, d\mu^n = \int_Y \int_G F(g \x_0) \, d\a_\y\, d\nu.$$
\end{lem}

\pf
Let $\varphi:X^n\rightarrow X^n$ be defined by 
$\varphi(\x)=(x_1,x_1^{-1}x_2,\dots,x_1^{-1}x_n)$ and $\pi$ be the image of $\mu^n$ with respect to $\varphi$, i.e, $\pi(A)=\mu^n\{\varphi\in A\}$ for all Borel $A\subset X^n$. Then
\[
\int\widetilde{F}\left(\varphi(\x)\right)d\mu^n(\x)=\int\widetilde{F}(\y)d\pi(\y)
\]
for nonnegative Borel functions $\widetilde{F}$ on $X^n$. Taking $\widetilde{F}(\y)=F(y_1,y_1y_2,\dots,y_1y_n)$ yields
\begin{eqnarray*}
\int F(\x)d\mu^n(\x)&=&\int F(y_1,y_1y_2,\dots,y_1y_n)d\pi(\y)\\
&=&\int d\widetilde{\nu}(y_2,\dots,y_n)\int F(y_1,y_1y_2,\dots,y_1y_n)d\widetilde{\alpha}_{(y_2,\dots,y_n)}(y_1),
\end{eqnarray*}
for some measures $\widetilde{\nu}$ and $\widetilde{\alpha}$. The right-hand side can then be written as 
\[
\int_{X_0}d\nu(\x_0)\int_GF(g\x_0)d\alpha_{\x_0}(g),
\]
where $\nu$ is an image of $\widetilde{\nu}$ with respect to the function $(y_2,\dots,y_n)\mapsto(1,y_2,\dots,y_n)$ and $\alpha_{\x_0}=\widetilde{\alpha}_{(y_2,\dots,y_n)}$ for $\x_0=(1,y_2,\dots,y_n)$, completing the proof.
\qed

\begin{lem}\label{Lem:Compact}
If $s$ is a shift-invariant ($n$-sample) estimator then 
$$Q(s)=\int_Y \int_G \chi_{ B ( s(\x_0)^{-1} ) } \,d\a_{\x_0} \, d\nu.$$
\end{lem}

\pf
Since $s$ is shift-invariant, its quality can be computed at the identity.
Thus $Q(s)=Q^1(s)=\mu^n(s^{-1}(B(1)))=\int_{X^n} \chi_{s^{-1}(B(1))} \,d\mu^n$.
By Lemma \ref{Lem:Decompose}, this integral can be decomposed as
$$\int_Y \int_G \chi_{s^{-1}(B(1))}(g \x_0) \,d\a_{\x_0} \,d\nu.$$
Now, note that $g \x_0\in s^{-1}(B(1))$ if and only if 
$gs(\x_0)\in B(1)$ if and only if $g\in B(s(\x_0)^{-1})$.
Thus the integral above is the same as the one in the statement
of the lemma, and we are done.
\qed

We are now ready to prove the result.  Note that we
do not prove that optimal estimators exist---only that if they exist, then one of
them is shift-invariant.

\begin{thm}\label{Thm:Compact}
Let $G=X$ be a compact group, let $\delta$ and $n$ be given, and let $\mu$ 
be defined by a density function.  If $e:X^n\to G$ is any estimator then there 
exists a shift-invariant estimator $s$ with $Q(s)\geq Q(e)$.
\end{thm}

\pf
Let $e:X^n\to G$ be any estimator.  For each group element $\gamma\in G$,
we define a shift-invariant estimator $s_\g$ that agrees with $e$ on the coset
$\gamma X_0$:
$$s_\g(g,gx_2,\ldots,gx_n)=g\g^{-1}e(\g,\g x_2,\ldots,\g x_n).$$
We will show that there exists $\g$ such that $Q(s_\g)\geq Q(e)$.

Denote by $\rho$ the invariant
(Haar) measure on $G$.  Since $Q(e)$ is defined as an infimum, we have
\begin{eqnarray}\label{Eqn:Triple}
Q(e)\leq \int_{\th\in G} Q^\th(e) \,d\rho
&=&\int_{\th\in G} \int_{X^n} \chi_{\th^{-1}e^{-1}(B(\th))}(\x)\,d\mu^n \,d\rho \nonumber \\
&=& \int_{X^n} \int_{\th\in G}\chi_{\th^{-1}e^{-1}(B(\th))}(\x)\,d\rho \,d\mu^n \nonumber \\
&=& \int_Y \int_G \int_{\th\in G}\chi_{\th^{-1}e^{-1}(B(\th))}(g\x_0)\,d\rho \,d\a_{\x_0} \,d\nu,
\end{eqnarray}
where the last equality comes from Lemma \ref{Lem:Decompose}.
The condition that $g \x_0\in\th^{-1}e^{-1}(B(\th))$ is equivalent to $d(e(\th g\x_0),\th)<\d$.
Now we make the substitution $\gamma=\th g$.  Thus $\th=\gamma g^{-1}$,
and the condition becomes $d(e(\gamma \x_0),\gamma g^{-1})<\d$,
or, by invariance of the metric, $d(\gamma^{-1} e(\gamma \x_0),g^{-1})<\d$.
This says that $g^{-1}\in B(\gamma^{-1}e(\gamma \x_0))$, or
equivalently, $g\in B(e(\gamma \x_0)^{-1} \gamma)$.

This allows us to rewrite the triple integral (\ref{Eqn:Triple}), using the
measure-preserving transformation $\th\mapsto \g=\th g$, as
\begin{eqnarray*}
\int_Y \int_G \int_{\th\in G}\chi_{\th^{-1}e^{-1}(B(\th))}
\,d\rho\, d\a_{\x_0}\, d\nu &=&
\int_Y \int_G \int_{\gamma\in G}
\chi_{B(e(\gamma \x_0)^{-1}\gamma)} \,d\rho \,d\a_{\x_0} \,d\nu \\
&=& 
\int_{\gamma\in G} \left( \int_Y \int_G 
\chi_{B(e(\gamma \x_0)^{-1}\gamma)} \,d\a_{\x_0} \,d\nu \right) d\rho \\
\end{eqnarray*}
Now, comparing with Lemma \ref{Lem:Compact}, we see that the inner integral
above is exactly the quality of the shift-invariant estimator $s_\g$.

We therefore have $$Q(e)\leq \int_\gamma Q(s_\g) \,d\rho,$$
or in other words, the average quality of the shift-invariant estimators $\{s_\g\}$ is at least
$Q(e)$.  Therefore at least one of the $s_\g$ satisfies $Q(s_\g)\geq Q(e)$.
\qed


\section{A non-shift-invariant example}
\label{Sec:Counterexample}


In the following example, suggested by L.~Schulman, the optimal shift-invariant
estimator is not optimal.  This provides an interesting
complement to Conjecture 1. 
Lehman and Casella \cite[Section 5.3]{lehman} also give examples of this phenomenon.

Let $X$ be the infinite trivalent tree, which we view as the Cayley graph
of the group $G=(\Z/2\Z) \ast (\Z/2\Z) \ast (\Z/2\Z) = 
\langle a,b,c \mid a^2=b^2=c^2=1 \rangle$.
In words, $G$ is a discrete group generated by three
elements $a,b,$ and $c$, each of order two and with no other relations.
Each non-identity element of $G$ can be written uniquely as
a finite word in the letters $a,b,c$ with no letter appearing twice in
a row; we refer to such a word as the \emph{reduced form} of
the group element.  (We write 1 for the identity element of $G$.)  
Multiplication in the group is performed by concatenating
words and then canceling any repeated letters in pairs.  Evidently
$G$ is infinite.  The Cayley graph $X$ is a graph with one vertex
labeled by each group element of $G$, and with an edge 
joining vertices $w$ and $x$ if and only if 
$w=xa$, $w=xb$, or $w=xc$.  Note that this relation is symmetric,
since $a$, $b$, and $c$ each have order 2.  Each vertex of $X$
has valence 3, and $X$ is connected and contains no circuits, 
i.e.~it is a tree.  Finally, $X$ becomes a metric space by declaring
each edge to have length one.

Because of how we defined the edges
of $X$, $G$ acts naturally on the left of $X$:  given $g\in G$, the map 
$g:X\to X$ defined by $g(x)=gx$ is an isometry of $X$.
So if $\d>0$ is given, $\mu$ is a probability distribution, $\th\in G$, 
and $e$ is an estimator, then the shift $\mu_\th$ and the quality $Q_\d(e)$
are defined as usual by (\ref{Def:Shift}) and (\ref{Def:Quality}). 

We are ready to present the example.
Suppose $0<d<1$ is fixed, and let $\mu$ be the probability 
distribution with atoms of weight $1/3$ at the three vertices $a,b,c$.
Thus for $\th\in G$, the distribution $\mu_\th$ has atoms of weight 
$1/3$ at the three neighbors $\th a,\th b,\th c$ of the vertex $\th$ in $X$.

\begin{example}
There is an optimal one-sample estimator with quality $2/3$, but the 
quality of any shift-invariant one-sample estimator is at most $1/3$.
\end{example}

\pf
Consider the one-sample estimator $e$ that truncates the last letter
of the sample (unless the sample is the identity, in which case we 
arbitrarily assign the value $a$).  That is, for a vertex $x$ of $X$,
\begin{equation*}
e(x)=
\begin{cases}
w \text{ if } x=w\ell \text{ is reduced, and } \ell=a,b, \text{ or } c \\
a \text{ if } x=1.
\end{cases}
\end{equation*}
Geometrically, this estimator takes a sample $x$ and, unless $x=1$, 
guesses that the shift is the (unique) neighbor of $x$ that is closer to 1.

We compute the quality of $e$.  Note $Q^1(e)=1$, because if $\th=1$ 
then the sample will be $a,b,$ or $c$, and the estimator is guaranteed 
to guess correctly.  In fact $Q^a(e)=1$ also, as is easily verified. 
For any other shift $\th$, the sample is $\th \ell$ for 
$\ell=a,b,$ or $c$, and the estimator guesses correctly exactly when
$\ell$ differs from the last letter of $\th$.  So $Q^\th(e)=2/3$, and 
$Q(e)=\inf_\th Q^\th(e)=2/3$.

It is easy to see that this estimator is optimal.  Suppose $e'$ is 
another estimator and $Q(e')>2/3$.  Since each local quality
$Q^\th(e')$ is either 0, 1/3, 2/3, or 1, we must have $Q^\th(e')=1$
for all $\th$.  This means $e'$ always guesses correctly.  But 
since there are different values of $\th$ that can produce
the same sample, this is impossible.

Observe that the estimator $e$ above is neither left- nor right-invariant.  
For instance right-invariance fails, as $e(ba\cdot a) = e(b) =1 \neq ba = e(ba)\cdot a$,
and the same example shows the failure of left-invariance:
$e(b\cdot aa) = id \neq  ba = be(aa)$.

Indeed, we conclude by showing that the quality of any shift-invariant
one-sample estimator $e'$ is at most 1/3.  Suppose $e'(1)=w$.  If $e'$ 
is left-invariant, it follows that $e'(x)=xw$ for all $x$; if $e'$ is right-invariant 
it follows that $e'(x)=wx$ for all $x$.  

Since $\d<1$, the quality of $e'$ at $\th=1$ is equal to the probability 
that $e'(x)=1$, given that $x$ was sampled from $\mu$.  With equal probability $x$
is $a$, $b$, or $c$; since at most one of $wa$, $wb$, $wc$ and one
of $aw$, $bw$, $cw$ can equal 1, we conclude that $Q(e')\leq Q^1(e')\leq 1/3$.
\qed

We remark that this example readily generalizes
to other finitely generated groups with infinitely many ends:  the key
is that $e$ is a two-to-one map but with only one sample, 
a shift-invariant estimator is necessarily one-to-one.


\end{document}